\documentclass[12pt]{article}

\usepackage{amssymb}
\usepackage{amsmath}
\usepackage{graphicx}
\usepackage{hyperref}

\numberwithin{equation}{section}

\begin{document}

\newcommand{\bd}{\begin{displaymath}}
\newcommand{\ed}{\end{displaymath}}
\newcommand{\ds}{\displaystyle}
\newcommand{\bp}{\underline{\bf Proof}:\ }
\newcommand{\ep}{{\hfill $\Box$}\\ }
\newcommand{\be}{\begin{equation}}
\newcommand{\ee}{\end{equation}}
\newcommand{\ba}{\begin{array}}
\newcommand{\ea}{\end{array}}
\newcommand{\bea}{\begin{eqnarray}}
\newcommand{\eea}{\end{eqnarray}}
\newcommand{\nt}{\noindent}

\newtheorem{0}{DEFINITION}[section]
\newtheorem{1}{LEMMA}[section]
\newtheorem{2}{THEOREM}[section]
\newtheorem{3}{COROLLARY}[section]
\newtheorem{4}{PROPOSITION}[section]
\newtheorem{5}{REMARK}[section]
\newtheorem{6}{EXAMPLE}[section]
\newtheorem{7}{ALGORITHM}[section]
\newtheorem{8}{CONJECTURE}[section]

\title{On Limiting Probability Distributions of Higher Order Markov Chains}
\author{
Lixing Han\thanks{Department of Mathematics, University of Michigan-Flint, Flint, MI 48502, USA. Email: \texttt{lxhan@umich.edu}}
\and
Jianhong Xu\thanks{School of Mathematical and Statistical Sciences, Southern Illinois University Carbondale, Carbondale, IL 62901, USA. Email: \texttt{jxu@math.siu.edu}} 
}

\maketitle

\begin{abstract}
The limiting probability distribution is one of the key characteristics of a Markov chain since it shows its long-term behavior. In this paper, for a higher order Markov chain, we establish some properties related to its exact limiting probability distribution, including a sufficient condition for the existence of such a distribution. Our results extend the corresponding conclusions on first order chains. Besides, they complement the existing results concerning higher order chains which rely on approximation schemes or two-phase power iterations. Several illustrative example are also given.
\end{abstract}

\nt {\bf Keywords}: higher order Markov chains, limiting probability distributions, stationary probability distributions, transition tensors, regularity


\section{Introduction}
\label{intro}
\setcounter{equation}{0}

Let $m \ge 2$ be an integer. An $(m-1)$th order Markov chain is a stochastic process $X=\{X_t : t=1, 2, \ldots\}$, where each $X_t$ takes values in state space $S=\{1, 2, \ldots, n\}$, such that 
\be
\label{mar}
\ba{l}
\Pr(X_{t+1}=i_1 | X_t=i_2, \ldots, X_{t-m+2}=i_m, \ldots, X_1=i_{t+1})\\[.5em]
=\Pr(X_{t+1}=i_1 | X_t=i_2, \ldots, X_{t-m+2}=i_m)
\ea
\ee
for all $t \ge m-1$ and $i_1, \ldots, i_m, \ldots, i_{t+1} \in S$. Throughout this paper, the chain $X$ is assumed to be homogeneous, meaning that the probability as in (\ref{mar}) is independent of $t$. Hence, we can denote 
$$\Pr(X_{t+1}=i_1 | X_t=i_2, \ldots, X_{t-m+2}=i_m)=p_{i_1i_2\ldots i_m},$$ 
for any $t \ge m-1$ and $i_1, i_2, \ldots, i_m \in S$. This is the transition probability from states $(i_2, \ldots, i_m)$ to state $i_1$. The resulting $m$th order, $n$ dimensional tensor ${\cal P}=[p_{i_1i_2\ldots i_m}]$ is called the transition tensor. Since $0 \le p_{i_1i_2\ldots i_m} \le 1$ for any $i_1, i_2, \ldots, i_m \in S$ and $\ds \sum_{i_1=1}^n p_{i_1i_2\ldots i_m}=1$ for any $i_2, \ldots, i_m \in S$, $\cal P$ is called stochastic. The chain is said to be higher order whenever $m \ge 3$. In the first order case, the transition tensor $\cal P$ is often denoted by $P \in \mathbb R^{n \times n}$, which is called the transition matrix .

For the past decade or so, higher order Markov chains have drawn much attention in the literature, see \cite{BozHaj16, ChaZha13, CulPeaZha17, DinNgWei18, Gei17, GleLimYu15, HanWanXu22, HanXu24a, HanXu24b, HuQi14, HuaQi20, LiZha16, LiNg14, WuChu17}. A focal point in a majority of these existing works is the limiting probability distribution of a higher order chain. The following definitions, however, include the first order case too since it is also relevant here.

\begin{0}
\label{dist}
For an $(m-1)$th order, where $m \ge 2$, Markov chain $X=\{X_t : t=1, 2, \ldots\}$ on state space $S=\{1, 2, \ldots, n\}$,
\be
\label{dist_for}
x_t=[\Pr(X_t=1) ~\Pr(X_t=2) ~\ldots ~\Pr(X_t=n)]^T \in \mathbb R^n
\ee
is called the probability distribution of $X$ at time $t$.
\end{0}

\begin{0}
\label{lim_dist}
Under the assumptions of Definition \ref{dist}, if $\ds \lim_{t \rightarrow \infty}x_t$ exists and if this limit does not depend on the initial probability distributions, i.e., $x_1, x_2, \ldots, x_{m-1}$, then 
$$\pi=\lim_{t \rightarrow \infty}x_t$$ 
is called the limiting probability distribution of the chain $X$.
\end{0}

Clearly, $\pi$ satisfies $\pi \ge 0$ and $\ds \sum_{i=1}^n \pi_i=1$ and thus is itself a probability distribution. It is an important characteristic of the chain since it provides useful information regarding its long-term behavior. The limit $\ds \lim_{t \rightarrow \infty} x_t$, however, may not exist and may depend on the initial probability distributions. An example for this shall be given later after we introduce some necessary background material. It is natural, therefore, to ask under what conditions the limiting probability distribution exists and what properties it must satisfy. 

For a first order chain with transition matrix $P$, a well-known sufficient condition regarding the existence of the limiting probability distribution is $P$ being regular. Denote\footnote{We may also denote this as $P^{(k)}$.} $P^k=[p^{(k)}_{ij}] \in \mathbb R^{n \times n}$, where $k=1, 2, \ldots$. This $p^{(k)}_{ij}$ represents the $k$-step transition probability from states $j$ to $i$, that is, 
\be
\label{kstep}
p^{(k)}_{ij}=\Pr(X_{t+k}=i | X_t=j)
\ee
for all $t \ge 1$. Then, the transition matrix $P$ (as well as the chain) is called regular if there exists some $k \ge 1$ such that $P^k > 0$. Under this condition, $\ds \lim_{t \rightarrow \infty}x_t=\pi$ exists, $\pi > 0$, and $\pi$ is independent of the initial probability distribution $x_1$. Moreover, $\ds \lim_{k \rightarrow \infty}P^k=\pi e^T$, where, and in what follows, $e$ is the vector of all ones whose dimension can be determined by the context. We refer the reader to \cite{Hun83, KemSne60} for the existing results on regular first order chains.

For the regular first order case, the limiting probability distribution is equivalent to the so-called stationary probability distribution. This appears to be the origin that the latter term has been used interchangeably in the literature on higher order chains. We, however, choose the former here because in the higher order case, the limiting probability distribution lacks an interpretation of stationarity as it comes to the evolution of the chain.

In light of the aforementioned results on first order chains, we now raise the question of whether parallel results can be established on higher order chains. As far as we can see, this question has been partially answered thus far by certain approximation schemes or a two-phase procedure under some stronger requirement, see a brief survey in Section \ref{back}. In this paper, we shall provide a more complete answer to this question. In particular, our results concern the exact limiting probability distribution without resorting to any approximation schemes. This is the main difference between our results and a vast majority of the relevant existing ones.

\section{Background Material}
\label{back}
\setcounter{equation}{0}

Let us set the stage here first for the presentation of our results in the next section.

Let $m \ge 3$. We consider in the sequel an $(m-1)$th order Markov chain $X=\{X_t : t=1, 2, \ldots\}$ on state space $S=\{1, 2, \ldots, n\}$ with an $m$th order, $n$ dimensional transition tensor ${\cal P}=[p_{i_1i_2\ldots i_m}]$.

For $k \ge 1$, similar to (\ref{kstep}), the $k$-step transition probability from states $(i_2,\ldots,i_m)$ to state $i_1$ is given by 
$$p^{(k)}_{i_1i_2\ldots i_m}=\Pr(X_{t+k}=i_1 | X_t=i_2, \ldots, X_{t-m+2}=i_m)$$
for any $t \ge m-1$. In particular, $p^{(1)}_{i_1\ldots i_m}=p_{i_1\ldots i_m}$. By convention, we define 
$$p^{(0)}_{i_1i_2\ldots i_m}=\delta_{i_1i_2\ldots i_m}=\left\{
\ba{cl}
1, & i_1=i_2;\\
0, & i_1 \ne i_2.
\ea\right.$$
These $k$-step transition probabilities form a stochastic tensor too, which is denoted by ${\cal P}^{(k)}=[p^{(k)}_{i_1\ldots i_m}]$. The special case of $k=0$ leads to the identity tensor ${\cal I}=[\delta_{i_1i_2\ldots i_m}]$ --- note that this identity tensor is different from those in \cite{KilMar11, QiLuo17}. In addition, it is known that the $k$-step transition probabilities can be determined by the following recurrence formula \cite{BozHaj16, HanWanXu22}:
\be
\label{pk}
p^{(k+1)}_{i_1i_2\ldots i_m}=\sum_{j=1}^n p^{(k)}_{i_1ji_2\ldots i_{m-1}}p_{ji_2\ldots i_m}, ~k=0, 1, 2, \ldots.
\ee
Using the symbol in \cite{HanXu24a}, (\ref{pk}) can be written as a tensor product ${\cal P}^{(k+1)}={\cal P}^{(k)} \boxtimes {\cal P}$. Especially, the identity tensor satisfies ${\cal P}={\cal I} \boxtimes {\cal P}$. In this sense, ${\cal P}^{(k)}$ can be thought of as the $k$th power of $\cal P$. The operation $\boxtimes$ is clearly distributive, yet it is neither commutative nor associative.

Next, we recall the notions of irreducibility, ergodicity, and regularity of the transition tensor $\cal P$ (as well as the higher order chain $X$). These can be extended to nonnegative tensors too, but we shall restrict ourselves to the stochastic scenario here. Specifically, $\cal P$ is said to be :

\begin{itemize}
\item{irreducible if for each $J$ such that $\emptyset \ne J \subsetneq S$, $p_{i_1i_2\ldots i_m}>0$ for some $i_1 \in J$ and some $i_2, \ldots, i_m \in J^c$.}
\item{ergodic if given any $i_1, i_2, \ldots, i_m \in S$, there exists $k \ge 1$, which may depend on $i_1, i_2, \ldots, i_m$, such that $p^{(k)}_{i_1i_2\ldots i_m}>0$.}
\item{regular if there exists some $k \ge 1$ such that $p^{(k)}_{i_1i_2\ldots i_m} > 0$ for every $i_1, i_2, \ldots, i_m \in S$.}
\end{itemize}
Obviously, regularity implies ergodicity. In the meantime, it is also known that ergodicity implies irreducibility but the reverse does not hold \cite{HanWanXu22}. The latter is one of the differences between the higher order and the first order cases because irreducibility and ergodicity are known to be equivalent on a first order chain.

Now, we review a systematical way of converting the higher order chain $X$ to a first order one. Denote $T=\{i_1i_2\ldots i_{m-1} : i_1, i_2, \ldots, i_{m-1} \in S\}$. It is the collection of all multi-indices of length $m-1$. Such multi-indices are arranged via linear indexing \cite{MarShaLar13}. When $m=4$ and $n=2$, for example, 
$$T=\{111, 211, 121, 221, 112, 212, 122, 222\}.$$
Clearly, the cardinality of $T$ is $N=n^{m-1}$. Consider random variables $Y_t=(X_t, X_{t-1}, \ldots, X_{t-m+2})$ that take values in $T$. Specifically, we define $Y_t=i_1i_2\ldots i_{m-1}$ whenever 
$X_t=i_1, X_{t-1}=i_2, \ldots, X_{t-m+2}=i_{m-1}$.
Then, $Y=\{Y_t : t=m-1, m, \ldots\}$ is a first order chain on state space $T$ \cite{Hun83, Ios07}. It is also referred to as the reduced first order chain corresponding to $X$. The transition matrix $Q \in \mathbb R^{N \times N}$ of $Y$ is given by 
$$q_{i_1i_2\ldots i_{m-1}, j_2j_3\ldots j_m}=\Pr(Y_{t+1}=i_1i_2\ldots i_{m-1} | Y_t=j_2j_3\ldots j_m).$$
Observe that $q_{i_1i_2\ldots i_{m-1}, j_2j_3\ldots j_m}=p_{i_1i_2\ldots i_{m-1}j_m}$ whenever $i_\ell=j_\ell$ for all $\ell=2, 3, \ldots, m-1$; otherwise, $q_{i_1i_2\ldots i_{m-1}, j_2j_3\ldots j_m}=0$. Again, considering $m=4$ and $n=2$, we arrive at 
$$Q=\left[\ba{cccccccc}
p_{1111} & 0 & 0 & 0 & p_{1112} & 0 & 0 & 0\\
p_{2111} & 0 & 0 & 0 & p_{2112} & 0 & 0 & 0\\
0 & p_{1211} & 0 & 0 & 0 & p_{1212} & 0 & 0\\
0 & p_{2211} & 0 & 0 & 0 & p_{2212} & 0 & 0\\
0 & 0 & p_{1121} & 0 & 0 & 0 & p_{1122} & 0\\
0 & 0 & p_{2121} & 0 & 0 & 0 & p_{2122} & 0\\
0 & 0 & 0 & p_{1221} & 0 & 0 & 0 & p_{1222}\\
0 & 0 & 0 & p_{2221} & 0 & 0 & 0 & p_{2222}
\ea\right].$$
In $Q$, $p_{i_1i_2\ldots i_m}$ is on row $i_1+n(i_2-1)+\ldots +n^{m-2}(i_{m-1}-1)$
and column $i_2+n(i_3-1)+\ldots +n^{m-2}(i_m-1)$.
As we shall see soon, there is an easier way of forming $Q$.

In the same manner as (\ref{kstep}), the $k$-step transition probability of $Y$ from states $j_2j_3\ldots j_m$ to $i_1i_2\ldots i_{m-1}$ is given by 
\be
\label{kstep2}
\ba{rcl}
q^{(k)}_{i_1i_2\ldots i_{m-1}, j_2j_3\ldots j_m} & = & \Pr(X_{t+k}=i_1, X_{t+k-1}=i_2, \ldots, X_{t+k-m+2}=i_{m-1}\\[.8em]
 & & \hspace{3em} | X_t=j_2, X_{t-1}=j_3, \ldots, X_{t-m+2}=j_m).
 \ea
 \ee
 Clearly, for $1 \le k \le m-2$, $q^{(k)}_{i_1i_2\ldots i_{m-1}, j_2j_3 \ldots j_m}=0$ unless $i_{k+1}=j_2, i_{k+2}=j_3, \ldots, i_{m-1}=j_{m-k}$.

It should be mentioned that problems regarding the higher order chain $X$ cannot always be fully addressed by resorting to the reduced first order chain $Y$. To explain this, we quote several examples here: (1) $Y$ may be non-ergodic even though $X$ is ergodic \cite{HanXu24a}, (2) the mean first passage times of $X$ may not be obtained in a direct fashion from those of $Y$ \cite{HanXu24a}, and (3) a recurrent state in $X$ may no longer be a part of any recurrent state in $Y$ \cite{Ios07}. The reduced first order chain, nevertheless, can still be a quite useful tool when dealing with the limiting probability distribution problem.

In the sequel, we denote\footnote{In the first order case, this notation is consistent with that of the transition matrix.} by $P \in \mathbb R^{n \times N}$ the mode-$1$ matricization \cite{KolBad09} of the transition tensor $\cal P$. Simply put, $P$ consists of the frontal slices of $\cal P$ which are arranged side by side and ordered according to linear indexing. Taking the case $m=4$ and $n=2$, $P$ has the form 
$$P=\left[\ba{cc|cc|cc|cc}
p_{1111} & p_{1211} & p_{1121} & p_{1221} & p_{1112} & p_{1212} & p_{1122} & p_{1222}\\
p_{2111} & p_{2211} & p_{2121} & p_{2221} & p_{2112} & p_{2212} & p_{2122} & p_{2222}
\ea\right].$$
We can also refer to the entries of $P$ in multi-index form as $p_{i_1, i_2\ldots i_m}$. In $P$, $p_{i_1i_2\ldots i_m}$ is on row $i_1$ and column $i_2+n(i_3-1)+\ldots +n^{m-2}(i_m-1)$. In a similar vein, $P^{(k)}$ denotes the mode-$1$ matricization of ${\cal P}^{(k)}$. By convention, $P^{(0)}$ is the mode-$1$ matricization of the identity tensor $\cal I$, that is, 
\be
\label{p0}
P^{(0)}=[\underbrace{I_n ~I_n ~\ldots ~I_n}_{n^{m-2}}],
\ee
where $I_n \in \mathbb R^{n \times n}$ is the identity matrix.

Incidentally, the transition matrix $Q$ of the reduced first order chain is constructed in \cite{WuChu17} using the product of a block diagonal matrix and a permutation matrix. By the mode-$1$ matricization of $\cal P$, the formulation of $Q$ can be simplified as 
\be
\label{q}
Q=G \ast P,
\ee
where $$G=[\underbrace{I_{n^{m-2}} ~I_{n^{m-2}} ~\ldots ~I_{n^{m-2}}}_n]$$ and $\ast$ stands for the columnwise Khatri-Rao product \cite{SmiBroGel04}.

More importantly, the mode-$1$ matricization of $\cal P$ can bridge the probability distribution of the higher order chain $X$ to that of the reduced first order chain $Y$. Specifically, let $y_t=x_{t,t-1,\ldots,t-m+2} \in \mathbb R^N$ be the probability distribution of $Y$ at time $t$, whose $i_1i_2\ldots i_{m-1}$th entry is given by 
$$\Pr(Y_t=i_1i_2\ldots i_{m-1})=\Pr(X_t=i_1, X_{t-1}=i_2, \ldots, X_{t-m+2}=i_{m-1})$$
for any $i_1i_2\ldots i_{m-1} \in T$. As usual, these entries are ordered by linear indexing. When $m=4$ and $n=2$, for example, we obtain 
$$y_t=\left[\ba{c}
\Pr(X_t=1, X_{t-1}=1, X_{t-2}=1)\\
\Pr(X_t=2, X_{t-1}=1, X_{t-2}=1)\\
\Pr(X_t=1, X_{t-1}=2, X_{t-2}=1)\\
\Pr(X_t=2, X_{t-1}=2, X_{t-2}=1)\\
\Pr(X_t=1, X_{t-1}=1, X_{t-2}=2)\\
\Pr(X_t=2, X_{t-1}=1, X_{t-2}=2)\\
\Pr(X_t=1, X_{t-1}=2, X_{t-2}=2)\\
\Pr(X_t=2, X_{t-1}=2, X_{t-2}=2)
\ea\right].$$
By $x_t$ in (\ref{dist_for}) and $y_t$ above, it can be shown \cite{LiNg14} that 
\be
\label{dist_rel}
x_{t+1}=Py_t.
\ee
Formula (\ref{dist_rel}) is exactly what the two-phase procedure in \cite{WuChu17} is based on.

Now, with the foregoing background material at hand, we are ready to present a brief survey of the existing lines of inquiry regarding the limiting probability distribution of the higher order chain $X$.

The first (and the mainstream) approach assumes that for the reduced first order chain $Y$, if $\ds \lim_{t \rightarrow \infty} y_t$ exists, then such a limit can be approximated by  
\be
\label{appr}
\underbrace{z \otimes_K \cdots \otimes_K z}_{m-1}
\ee
for some $z \in \mathbb R^n$, where $\otimes_K$ stands for the Kronecker product. This approach is adopted in most existing works on higher order chains; see, for example, \cite{BozHaj16, ChaZha13, CulPeaZha17, DinNgWei18, FasTud20, GauTudHei19, GleLimYu15, HuQi14, HuaQi20, LiZha16, LiNg14}. It leads to the interesting $Z$-eigenvector problem. It is known \cite{LiNg14} that when $\cal P$ is irreducible, $z$ can be determined as some normalized positive $Z$-eigenvector. As pointed out in \cite{WuChu17}, however, the approximation in (\ref{appr}) may not capture the true nature of $\ds \lim_{t \rightarrow \infty} y_t$.

The second approach \cite{WuChu17} involves a two-phase power iteration scheme based on the connection between the higher order chain $X$ and its reduced first order chain $Y$. It assumes that $Q$ is regular, and so $\lambda=1$ is the only dominant eigenvalue of $Q$ on the unit circle and is simple. It follows that $\ds \lim_{t \rightarrow \infty}y_t$ is the eigenvector, normalized using $\| \cdot \|_1$, of $Q$ associated with $\lambda$. Thus, the limiting probability distribution of $X$ can be determined by, for example, (\ref{dist_rel}). As we shall demonstrate, however, the limiting probability distribution of $X$ may still exist despite $Q$ being non-regular.

The third approach \cite{Gei17} deals with the higher order chain $X$ by a certain approximation scheme to its transition probabilities. The development also hinges on a direct extension of the notion of recurrent states from the first order case to the higher order case. Such an extension has been shown  \cite{HanXu24b} to be an issue. Incidentally, this issue manifests another difference between the first order and the higher order cases.

Motivated by the above works we have reviewed, especially by \cite{WuChu17}, our main goals here are to further investigate the properties of the actual (not approximated) limiting probability distribution of the higher order chain $X$ and to establish a new sufficient condition that guarantees the existence of such a limiting probability distribution.

\section{Main Results}
\label{main}
\setcounter{equation}{0}

With the background material we have presented, let us begin this section with one example to illustrate that, in general, $\ds \lim_{t \rightarrow \infty} x_t$ may depend on the initial probability distributions and, therefore, the limiting probability distribution may not exist.
\begin{6}
\label{ex1}
{\rm Take a second order chain $X$ on state space $\{1, 2\}$ whose transition tensor $\cal P$ is as follows:
$${\cal P}(:,:,1)=\left[\ba{cc}
1 & 1/2\\
0 & 1/2
\ea
\right], ~{\cal P}(:,:,2)=\left[\ba{cc}
1/2 & 0\\
1/2 & 1
\ea
\right].$$
If we start with $X_1=1$ and $X_2=1$, then $y_t=[1 ~0 ~0 ~0]^T$ for any $t \ge 1$ and so $\ds \lim_{t \rightarrow \infty} x_t=[1 ~0]^T$ by (\ref{dist_rel}). On the other hand, if we start with $X_1=2$ and $X_2=2$, then $y_t=[0 ~0 ~0 ~1]^T$ for any $t \ge 1$ and so $\ds \lim_{t \rightarrow \infty} x_t=[0 ~1]^T$.}
\end{6}

In the rest of this paper, we assume $m \ge 3$ unless otherwise stated.

For a second order chain, the following technical lemma can be seen in \cite{Vla85}. We now extend it to the more general higher order chain $X$. It should be pointed out, however, that such an extension from the second order case to the general higher order case is usually not straightforward, as evidenced by a shortage of results for the latter. A key reason for this situation lies in the fact that the general case is often far less tractable.

\begin{1}
\label{komo}
Let $X$ be an $(m-1)$th order Markov chain on state space $S$ whose  transition tensor is $\cal P$. Assume $u, v = 0, 1, 2, \ldots$. Let $P^{(u)}$ and $Q$ be the mode-$1$ matricization of ${\cal P}^{(u)}$ and the transition matrix of the reduced first order chain $Y$, respectively. Then, 
\be
\label{komo_for}
P^{(u+v)}=P^{(u)}Q^v.
\ee
In particular, we have 
\be
\label{komo_for0}
P^{(v)}=P^{(0)}Q^v
\ee
for any $v=0, 1, 2, \ldots$.
\end{1}
\bp
Clearly, (\ref{komo_for}) holds when $u \ge 0$ and $v=0$. Assume now $u=0$ but $v \ge 1$. Note that on row $i_1$, $P^{(0)}$ has ones at, in multi-index form, columns $i_1j_2\ldots j_{m-1}$ for any $j_2, \ldots, j_{m-1} \in S$ and zeros elsewhere.
Hence, 
$$
\ba{rcl}
(P^{(0)}Q^v)_{i_1,i_2\ldots i_m} & = &  \displaystyle \sum_{j_1\ldots j_{m-1} \in T} p^{(0)}_{i_1,j_1\ldots j_{m-1}}q^{(v)}_{j_1\ldots j_{m-1}, i_2\ldots i_m}\\[1.6em]
 & = & \displaystyle \sum_{j_2, \ldots, j_{m-1} \in S} q^{(v)}_{i_1j_2\ldots j_{m-1}, i_2\ldots i_m}\\[1.6em]
& = & \displaystyle \sum_{j_2, \ldots, j_{m-1} \in S}\Pr(X_{t+v}=i_1, X_{t+v-1}=j_2, \ldots, X_{t+v-m+2}=j_{m-1}\\
 & & \hspace{6em} | X_t=i_2, X_{t-1}=i_3, \ldots, X_{t-m+2}=i_m)\\[.8em]
 & = & p^{(v)}_{i_1, i_2\ldots i_m},
\ea$$
where, according to the remark below (\ref{kstep2}), the third summation above is effectively over $j_2\ldots j_v$ only when $2 \le v \le m-2$ and is equal to $p_{i_1,i_2\ldots i_m}$ directly when $v=1$.

Similarly, for $u, v \ge 1$, we have 
$$\ba{l}
(P^{(u)}Q^v)_{i_1, i_2\ldots i_m}\\[.8em]
 = \displaystyle \sum_{j_1\ldots j_{m-1} \in T}p^{(u)}_{i_1,j_1\ldots j_{m-1}}q^{(v)}_{j_1\ldots j_{m-1}, i_2\ldots i_m}\\[1.6em]
 = \displaystyle \sum_{j_1\ldots j_{m-1} \in T}\Pr(X_{t+u+v}=i_1 | X_{t+v}=j_1, \ldots, X_{t+v-m+2}=j_{m-1})\\
 \hspace{4em} \cdot \Pr(X_{t+v}=j_1, \ldots, X_{t+v-m+2}=j_{m-1} | X_t=i_2, \ldots, X_{t-m+2}=i_m)\\[.8em]
 = \displaystyle \sum_{j_1\ldots j_{m-1} \in T}\Pr(X_{t+u+v}=i_1 | X_{t+v}=j_1, \ldots, X_{t+v-m+2}=j_{m-1},\\[.8em]
 \hspace{12em} X_t=i_2, \ldots, X_{t-m+2}=i_m)\\[.8em]
 \hspace{4em} \cdot \Pr(X_{t+v}=j_1, \ldots, X_{t+v-m+2}=j_{m-1} | X_t=i_2, \ldots, X_{t-m+2}=i_m)\\[.8em] 
 = \displaystyle \sum_{j_1\ldots j_{m-1} \in T}\Pr(X_{t+u+v}=i_1, X_{t+v}=j_1, \ldots, X_{t+v-m+2}=j_{m-1}\\[.8em]
 \hspace{7em} | X_t=i_2, \ldots, X_{t-m+2}=i_m)\\[1em]
 = p^{(u+v)}_{i_1i_2\ldots i_m},
\ea$$
where in the third equality, the higher order Markovian property in (\ref{mar}) is used, while in the fourth equality, the remark below (\ref{kstep2}) should again be noted for the case $1 \le v \le m-2$.
\ep

The above lemma shows a connection between the power of $Q$ and the power of $\cal P$ in matricized form. It has several applications, including in the proof of Theorem \ref{conn} later. Here, we first prove a sufficient condition for $\cal P$ to be regular, which generalizes a similar result in \cite{Vla85} for the second order case:

\begin{2}
\label{test}
Let $\cal P$ be the transition tensor of the higher order chain $X$ and $Q$ be the transition matrix of its reduced first order chain $Y$. If $Q$ is regular, then so is $\cal P$.
\end{2}
\bp
Since $Q$ is regular, there exists $k \ge 1$ such that $Q^k>0$. The conclusion follows from (\ref{komo_for0}).
\ep

Theorem \ref{test} provides a way of determining the regularity of $\cal P$ through $Q$, which allows us to make use of the available results in this regard \cite{HorJoh85}. In general, the converse of Theorem \ref{test} is not true. We, however, have the following for the case ${\cal P}>0$:

\begin{2}
For the $(m-1)$th order Markov chain $X$ on state space $S$ with transition tensor ${\cal P}=[p_{i_1i_2\ldots i_m}]>0$, denote $Q^k=[q^{(k)}_{i_1i_2\ldots i_{m-1},j_1j_2\ldots j_{m-1}}]$ using multi-index form. Then $$q^{(m-1)}_{i_1i_2\ldots i_{m-1},j_1j_2\ldots j_{m-1}}=\underbrace{p_{i_1i_2\ldots i_{m-1}j_1}p_{i_2i_3\ldots i_{m-1}j_1j_2} \cdots p_{i_{m-1}j_1j_2\ldots j_{m-1}}}_{m-1}>0,$$ which implies that $Q$ is regular. 
\end{2}
\bp
We claim that for any $2 \le \ell \le m-1$, 
\be
\label{ql}
\begin{array}{l}
q^{(\ell)}_{i_1i_2\ldots i_{m-1},j_1j_2\ldots j_{m-1}}\\
=\left\{\begin{array}{cl}
\underbrace{p_{i_1i_2\ldots i_{m-1}j_{m-\ell}}p_{i_2i_3\ldots i_{m-1}j_{m-\ell}j_{m-\ell+1}}\cdots p_{i_\ell i_{\ell+1}\ldots i_{m-1}j_{m-\ell}j_{m-\ell+1}\ldots j_{m-1}}}_\ell, & \begin{array}{l}
j_1=i_{\ell+1},\\
j_2=i_{\ell+2}, \ldots,\\
j_{m-\ell-1}=i_{m-1};
\end{array}\\
0, & ~{\rm otherwise}.
\end{array}\right.
\end{array}
\ee
Note that the number of conditions for having a positive $q^{(\ell)}_{i_1i_2\ldots i_{m-1},j_1j_2\ldots j_{m-1}}$ decreases as $\ell$ increases. Any such condition is interpreted here as null when any subscript is out of range. Let us verify the above by induction.

When $\ell=2$, 
$$\begin{array}{rcl}
q^{(2)}_{i_1i_2\ldots i_{m-1},j_1j_2\ldots j_{m-1}} & = & \displaystyle \sum_{k_1k_2\ldots k_{m-1} \in T}q_{i_1i_2\ldots i_{m-1},k_1k_2\ldots k_{m-1}}q_{k_1k_2\ldots k_{m-1},j_1j_2\ldots j_{m-1}}\\
 & = & q_{i_1i_2\ldots i_{m-1},i_2i_3\ldots i_{m-1}j_{m-2}}q_{i_2i_3\ldots i_{m-1}j_{m-2},j_1j_2\ldots j_{m-1}}
\end{array}$$
since $k_1=i_2, k_2=i_3, \ldots, k_{m-2}=i_{m-1}$, and $k_{m-1}=j_{m-2}$; otherwise, either $q_{i_1i_2\ldots i_{m-1},k_1k_2\ldots k_{m-1}}=0$ or $q_{k_1k_2\ldots k_{m-1},j_1j_2\ldots j_{m-1}}=0$. Thus, 
$$q^{(2)}_{i_1i_2\ldots i_{m-1},j_1j_2\ldots j_{m-1}}=\left\{\begin{array}{cl}
p_{i_1i_2\ldots i_{m-1}j_{m-2}}p_{i_2i_3\ldots i_{m-1}j_{m-2}j_{m-1}}, & j_1=i_3, j_2=i_4, \ldots, j_{m-3}=i_{m-1};\\
0, & {\rm otherwise}. 
\end{array}\right.$$
Suppose now that the claim holds for any $2 \le \ell \le m-2$. Then, 
$$\begin{array}{rcl}
q^{(\ell+1)}_{i_1i_2\ldots i_{m-1},j_1j_2\ldots j_{m-1}} & = & \displaystyle \sum_{k_1k_2\ldots k_{m-1} \in T}q_{i_1i_2\ldots i_{m-1},k_1k_2\ldots k_{m-1}}q^{(\ell)}_{k_1k_2\ldots k_{m-1},j_1j_2\ldots j_{m-1}}\\
 & = & p_{i_1i_2\ldots i_{m-1}j_{m-\ell-1}}q^{(\ell)}_{i_2i_3\ldots i_{m-1}j_{m-\ell-1},j_1j_2\ldots j_{m-1}},
\end{array}$$
where we note that $k_1=i_2, k_2=i_3, \ldots, k_{\ell+1}=i_{\ell+2}, \ldots, k_{m-2}=i_{m-1}$, and, using the inductive hypothesis, $k_{m-1}=j_{m-\ell-1}$. Furthermore, using the inductive hypothesis again, 
$$\begin{array}{l}
q^{(\ell+1)}_{i_1i_2\ldots i_{m-1},j_1j_2\ldots j_{m-1}}\\
=\left\{\begin{array}{cl}
p_{i_1i_2\ldots i_{m-1}j_{m-\ell-1}}\underbrace{p_{k_1k_2\ldots k_{m-1}j_{m-\ell}}\cdots p_{k_\ell k_{\ell+1}\ldots k_{m-1}j_{m-\ell}j_{m-\ell+1}\ldots j_{m-1}}}_{\ell}, & \begin{array}{l}
j_1=k_{\ell+1},\\
j_2=k_{\ell+2}, \ldots,\\
j_{m-\ell-1}=k_{m-1};
\end{array}\\
0, & ~{\rm otherwise},
\end{array}\right.\\

=\left\{\begin{array}{cl}
\underbrace{p_{i_1i_2\ldots i_{m-1}j_{m-\ell-1}}p_{i_2i_3\ldots i_{m-1}j_{m-\ell-1}j_{m-\ell}}\cdots p_{i_{\ell+1} i_{\ell+2}\ldots i_{m-1}j_{m-\ell-1}j_{m-\ell}\ldots j_{m-1}}}_{\ell+1}, & \begin{array}{l}
j_1=i_{\ell+2},\\
j_2=i_{\ell+3}, \ldots,\\
j_{m-\ell-2}=i_{m-1};
\end{array}\\
0, & ~{\rm otherwise}.
\end{array}\right.

\end{array}$$
The proof, therefore, is done.
\ep

In what follows, $\otimes$ stands for the outer product. We now state our main result as follows.

\begin{2}
\label{conv}
Consider the $(m-1)$th order Markov chain on state space $S$. Suppose that its transition tensor $\cal P$ is regular. Then, there exists $\pi \in \mathbb R^n$ satisfying $\pi_i > 0$ for each $i$ and $\ds \sum_{i=1}^n \pi_i = 1$ such that $\ds \lim_{k \rightarrow \infty} p^{(k)}_{ii_2\ldots i_m}=\pi_i$ for all $i_2, \ldots, i_m \in S$, i.e. 
\be
\label{lim_for}
\lim_{k \rightarrow \infty} {\cal P}^{(k)}=\pi \otimes \underbrace{e \otimes \cdots \otimes e}_{m-1}.
\ee
\end{2}
\bp
Suppose first that ${\cal P} > 0$. For any $x \in \mathbb R^n$, we define the product of $x$ with any mode-$1$ fiber of ${\cal P}^{(k)}$ as
\be
\label{xp}
x^T{\cal P}^{(k)}(:,i_2,\ldots,i_m)=\sum_{i=1}^n x_ip^{(k)}_{ii_2\ldots i_m}.
\ee
Denote $x^T{\cal P}^{(k)}(:,i_2,\ldots,i_m)$ as $y^{(k)}_{i_2\ldots i_m}$. Clearly, $y^{(k)}=[y^{(k)}_{i_2\ldots i_m}]$ is an $(m-1)$th order, $n$-dimensional tensor, but it can also be thought of as a vector whose entries are indexed by multi-indices $i_2\ldots i_m$. By (\ref{pk}), we see 
\be
\label{y}
y^{(k+1)}_{i_2\ldots i_m}=\sum_{j=1}^n p_{ji_2\ldots i_m}y^{(k)}_{ji_2\ldots i_{m-1}}.
\ee

Let $\epsilon=\min_{i_1\ldots i_m}p_{i_1\ldots i_m}$, $U_k=\max_{i_2\ldots i_m} y^{(k)}_{i_2\ldots i_m}$, and $L_k=\min_{i_2\ldots i_m} y^{(k)}_{i_2\ldots i_m}$. In particular, we set $U_0=\max_i x_i$ and $L_0=\min_i x_i$. Also, without loss of generality, we assume $0 < \epsilon \le 1/2$.
  
Let $\hat x \in \mathbb R^n$ be obtained from $x$ by replacing all $x_i$ with $U_0$ except one entry $x_\ell$ that equals $L_0$. Since $x \le \hat x$, we know that for all $i_2, \ldots, i_m$, $$y^{(1)}_{i_2\ldots i_m} \le \hat x^T{\cal P}(:,i_2,\ldots,i_m).$$
Concerning $\hat x^T{\cal P}(:,i_2,\ldots,i_m)$, we have 
\be
\begin{split}
\hat x^T{\cal P}(:,i_2,\ldots,i_m) & =[U_0e^T-(U_0-L_0)e_\ell^T]{\cal P}(:,i_2,\ldots,i_m)\\[.5em] \nonumber
 & =U_0-(U_0-L_0)p_{\ell i_2\ldots i_m}\\[.5em]
 & \le U_0-\epsilon (U_0-L_0),
\end{split}
\ee
which leads to $U_1 \le U_0-\epsilon (U_0-L_0)$. By a similar argument in which $\hat x$ is obtained from $x$ by replacing all $x_i$ with $L_0$ except one entry $x_\ell$ which is equal to $U_0$, we obtain $L_1 \ge L_0+\epsilon (U_0-L_0)$. Thus,
$$U_1-L_1 \le (1-2\epsilon)(U_0-L_0).$$

From (\ref{y}), we have $y^{(k+1)}_{i_2\ldots i_m} \le \sum_{j=1}^n p_{ji_2\ldots i_m} U_k = U_k$ and so $U_{k+1} \le U_k$, showing that $\{U_k\}$ is a non-increasing sequence. Similarly, $\{L_k\}$ is a non-decreasing sequence. Accordingly, $\{U_k - L_k\}$ is a non-increasing sequence. In order to prove the convergence of $\{U_k - L_k\}$, therefore, we only need to show that it has a convergent subsequence.

Applying (\ref{y}) repeatedly, we observe that for $k > m-2$,
$$y^{(k+1)}_{i_2\ldots i_m}=\sum_{j_{m-1}=1}^n p_{j_{m-1}i_2\ldots i_m}\sum_{j_{m-2}=1}^n p_{j_{m-2}j_{m-1}i_2\ldots i_{m-1}}\ldots \sum_{j_1=1}^n p_{j_1j_2\ldots j_{m-1}i_2}y^{(k-m+2)}_{j_1j_2\ldots j_{m-1}}.$$
Let $\hat y^{(k-m+2)}$ be the vector which is obtained from $y^{(k-m+2)}$ by replacing all its entries with $U_{k-m+2}$ except one entry $y^{(k-m+2)}_{\ell_1\ldots \ell_{m-1}}$ which is equal to $L_{k-m+2}$. Then,
\be
\begin{split}
y^{(k+1)}_{i_2\ldots i_m} & \le \sum_{j_{m-1}=1}^n p_{j_{m-1}i_2\ldots i_m}\sum_{j_{m-2}=1}^n p_{j_{m-2}j_{m-1}i_2\ldots i_{m-1}}\ldots \sum_{j_1=1}^n p_{j_1j_2\ldots j_{m-1}i_2}\hat y^{(k-m+2)}_{j_1j_2\ldots j_{m-1}}\\ \nonumber
& = \sum_{j_{m-1}=1}^n p_{j_{m-1}i_2\ldots i_m}\sum_{j_{m-2}=1}^n p_{j_{m-2}j_{m-1}i_2\ldots i_{m-1}}\ldots \sum_{j_1=1}^n p_{j_1j_2\ldots j_{m-1}i_2}U_{k-m+2}\\
& \hspace{.3in} -p_{\ell_{m-1}i_2\ldots i_m}p_{\ell_{m-2}\ell_{m-1}i_2\ldots i_{m-1}}\cdots p_{\ell_1\ell_2\ldots \ell_{m-1}i_2}(U_{k-m+2}-L_{k-m+2})\\
\\
& \le U_{k-m+2}-\epsilon^{m-1}(U_{k-m+2}-L_{k-m+2}),
  \end{split}
\ee
which leads to $$U_{k+1} \le U_{k-m+2}-\epsilon^{m-1}(U_{k-m+2}-L_{k-m+2}).$$
Similarly, we obtain 
$$L_{k+1} \ge L_{k-m+2}+\epsilon^{m-1}(U_{k-m+2}-L_{k-m+2}).$$
Together, these last two expressions yield 
$$U_{k+1}-L_{k+1} \le (1-2\epsilon^{m-1})(U_{k-m+2}-L_{k-m+2}).$$
In particular, on letting $k=m-1, 2m-2, 3m-3, \ldots$, we arrive at 
$$U_m-L_m \le (1-2\epsilon^{m-1})(U_1-L_1),$$
$$U_{2m-1}-L_{2m-1} \le (1-2\epsilon^{m-1})(U_m-L_m),$$
$$U_{3m-2}-L_{3m-2} \le (1-2\epsilon^{m-1})(U_{2m-1}-L_{2m-1}),$$
and so on. It follows that the subsequence $$U_0-L_0, U_1-L_1, U_m-L_m, U_{2m-1}-L_{2m-1}, U_{3m-2}-L_{3m-2}, \ldots$$ converges to zero, and hence so does the entire sequence $\{U_k-L_k\}$.

The above analysis also shows that  $$\{U_0, U_1,U_m,U_{2m-1},U_{3m-2},\ldots\}$$ is strictly decreasing while $$\{L_0, L_1,L_m,L_{2m-1},L_{3m-2},\ldots\}$$ is strictly increasing unless $U_k-L_k=0$ for some $k=0, 1, m, 2m-1, 3m-2, \ldots$. In particular, if we start with $U_0-L_0>0$, then there must be some $L_0 < c < U_0$ such that $\lim_{k \rightarrow \infty}(U_k-L_k)=0$, i.e. $\lim_{k \rightarrow \infty}y^{(k)}_{i_2\ldots i_m}=c$ for all $i_2, \ldots, i_m$.

Next, let us denote by $e_i$, $1 \le i \le n$, the $i$th column of an $n \times n$ identity matrix. For each $i$, we pick $x=e_i$ in (\ref{xp}) and apply the above convergence result. In this case, $L_0=0$ and $U_0=1$. It follows, therefore, that for each $i$, there exists $0 < \pi_i < 1$ such that $\lim_{k \rightarrow \infty}p^{(k)}_{ii_2\ldots i_m}=\pi_i$ for all $i_2, \ldots, i_m$. The fact that $\sum_{i=1}^n \pi_i=1$ is obvious because of the stochasticity of ${\cal P}^{(k)}$.

Finally, in the case that $\cal P$ is regular but ${\cal P} \not > 0$, there exists $K \ge 1$ such that ${\cal P}^{(K)}>0$. We set $\epsilon=\min_{i_1\ldots i_m}p^{(K)}_{i_1\ldots i_m}$ and consider the subsequence $$U_0-L_0, U_K-L_K, U_{mK}-L_{mK}, U_{(2m-1)K}-L_{(2m-1)K}, U_{(3m-2)K}-L_{(3m-2)K},\ldots.$$ The proof can be done in exactly the same fashion as the ${\cal P}>0$ case.
\ep

We comment that for the second order case, a result similar to Theorem \ref{conv} can be found in \cite{Vla85}. Our next result confirms that $\pi$ as in Theorem \ref{conv} is indeed the limiting probability distribution of the chain, which contains the long-term probabilities for the chain to be in various states. Moreover, these probabilities are independent of the initial probability distributions.

\begin{2}
\label{prob}
Let $X=\{X_t : t=1, 2, \ldots\}$ be an $(m-1)$th order Markov chain on state space $S$ whose transition tensor is $\cal P$. If there exists a nonnegative $\pi \in \mathbb R^n$ such that $\ds \sum_{i=1}^n \pi_i = 1$ and 
$\ds{\lim_{k \rightarrow \infty} {\cal P}^{(k)}=\pi \otimes \underbrace{e \otimes \cdots \otimes e}_{m-1}}$, then 
$$
\lim_{t \rightarrow \infty}\Pr(X_t=i)=\pi_i
$$
for each $i \in S$.
\end{2}
\bp
For any $k \ge 1$, we observe 
\be
\begin{split}
\Pr(X_{m+k-1}=i) & = \sum_{i_2, \ldots, i_m \in S} \Pr(X_{m+k-1}=i, X_{m-1}=i_2, \ldots, X_1=i_m)\\[.5em] \nonumber
 & = \sum_{i_2, \ldots, i_m \in S} \Pr(X_{m+k-1}=i | X_{m-1}=i_2, \ldots, X_1=i_m)\\
 & \hspace{.8in} \cdot \Pr(X_{m-1}=i_2, \ldots, X_1=i_m)\\[.5em]
 & = \sum_{i_2, \ldots, i_m \in S} p^{(k)}_{ii_2\ldots i_m}\Pr(X_{m-1}=i_2, \ldots, X_1=i_m).
\end{split}
\ee
Pushing $k \rightarrow \infty$ and noting $\ds \sum_{i_2, \ldots, i_m \in S} \Pr(X_{m-1}=i_2, \ldots, X_1=i_m)=1$, the proof is now complete.
\ep

We remark that Theorem \ref{prob} does not require $\pi$ to be strictly positive. This will allow us to explore in future works limiting behavior of the chain $X$ in other settings such as being a higher order absorbing chain.

We also point out that if there exists $y \in \mathbb R^N$ such that $\ds \lim_{k \rightarrow \infty}Q^k=y e^T$, then, by Lemma \ref{komo}, the mode-$1$ matricization of ${\cal P}^{(k)}$ satisfies $\ds \lim_{k \rightarrow \infty}P^{(k)}=(P^{(0)}y)e^T$. This translates exactly to $\ds{\lim_{k \rightarrow \infty} {\cal P}^{(k)}=\pi \otimes \underbrace{e \otimes \cdots \otimes e}_{m-1}}$, where $\pi=P^{(0)}y$ is the limiting probability distribution according to Theorem \ref{prob}. We shall further connect the limiting probability distribution to an eigenvector of $Q$ later in this section.

Along with Theorem \ref{prob}, Theorem \ref{conv} establishes a new sufficient condition, namely, $\cal P$ being regular, for the existence of the limiting probability distribution of the chain $X$. Moreover, it shows that in this case, the limiting probability distribution is positive. These extend the well-known results regarding first order chains as discussed in our introductory section.

Besides, we mention that without the condition of $\cal P$ being regular, the limiting probability distribution may not exist --- see the following example. Hence, such a sufficient condition is ``tight''. 
\begin{6}
{\rm
Let us consider a second order chain on the state space $\{1,2,3\}$ with transition tensor $${\cal P}(:,:,i_3)=\left[\ba{ccc}
0 & 1/2 & 0\\
1 & 0 & 1\\
0 & 1/2 & 0
\ea\right]$$
for any $i_3=1,2,3$. This $\cal P$ is ergodic but is non-regular. When $k$ is odd, we have ${\cal P}^{(k)}={\cal P}$. When $k$ is even, we have $${\cal P}^{(k)}(:,:,i_3)=\left[\ba{ccc}
1/2 & 0 & 1/2\\
0 & 1 & 0\\
1/2 & 0 & 1/2
\ea\right]$$
for $i_3=1, 2, 3$. As a result, $\ds \lim_{k \rightarrow \infty}{\cal P}^{(k)}$ does not exist.

To justify that this chain has no limiting probability distribution, let us assume $X_1=X_2=1$. Then, $x_t$ may be determined by (\ref{dist_rel}). Alternatively, we can proceed as follows.

By $p_{211}=1$, we see $X_3=2$, i.e. $x_3=[0 ~1 ~0]^T$. Next, using $p_{121}=p_{321}=1/2$, the chain will move to either state $1$ or state $3$ with equal probability, i.e. $x_4=[1/2 ~~0 ~~1/2]^T$. Continuing, since $p_{212}=p_{232}=1$, the chain will return to state $2$ in its next move, i.e. $x_5=[0 ~1 ~0]^T$. Hence, $x_t$ will bounce back and forth between $[0 ~1 ~0]^T$ and $[1/2 ~~0 ~~1/2]^T$ and will not converge as $t \rightarrow \infty$.}
\end{6}

The proof of Theorem \ref{prob} demonstrates how $\displaystyle{\lim_{k \rightarrow \infty} {\cal P}^{(k)}=\pi \otimes \underbrace{e \otimes \cdots \otimes e}_{m-1}}$ can ensure that the limiting probability distribution does not hinge on the initial probability distributions. The result below deals with the case when $\ds \lim_{k \rightarrow \infty}{\cal P}^{(k)}$ no longer has the special rank-one structure as in Theorem \ref{prob}.

\begin{2}
\label{prob2}
Consider an $(m-1)$th order Markov chain $X=\{X_t : t=1, 2, \ldots\}$ with state space $S$ and transition tensor $\cal P$. If there exists a tensor $\cal R$ such that $\ds \lim_{k \rightarrow \infty} {\cal P}^{(k)}=\cal R$, then for each $i \in S$, $\ds \lim_{t \rightarrow \infty}\Pr(X_t=i)$ exists. In general, however, such a limit depends on the initial probability distributions of $X$.
\end{2}
\bp
Similar to the proof of Theorem \ref{prob}, for any $k \ge 1$, we find 
$$\Pr(X_{m+k-1}=i)=\sum_{i_2,\ldots,i_m \in S} p^{(k)}_{ii_2\ldots i_m}\Pr(X_{m-1}=i_2, \ldots, X_1=i_m).$$
On letting ${\cal R}=[r_{i_1i_2\ldots i_m}]$, it follows that 
\be
\label{prob2_for}
\lim_{t \rightarrow \infty}\Pr(X_t=i) =\sum_{i_2, \ldots, i_m \in S} r_{ii_2\ldots i_m}\Pr(X_{m-1}=i_2, \ldots, X_1=i_m).
\ee
\ep

Now, in view of Theorem \ref{prob2}, we revisit Example \ref{ex1} to illustrate how $\ds \lim_{t \rightarrow \infty}\Pr(X_t=i)$ may depend on the initial probability distributions.
\begin{6}
{\rm 
Given the chain in Example \ref{ex1}, it is easy to verify that 
$$\lim_{k \rightarrow \infty}{\cal P}^{(k)}(:,:,1)=\left[\begin{array}{cc}
1 & 1/3\\
0 & 2/3
\end{array}\right] ~{\rm and} ~\lim_{k \rightarrow \infty}{\cal P}^{(k)}(:,:,2)=\left[\begin{array}{cc}
2/3 & 0\\
1/3 & 1
\end{array}\right].$$
According to (\ref{prob2_for}), we obtain that the entries of $\lim_{t \rightarrow \infty} x_t$ are given by 
$$\lim_{t \rightarrow \infty}\Pr(X_t=1)=\Pr(X_2=1,X_1=1)+\frac{\Pr(X_2=2,X_1=1)}{3}+\frac{2\Pr(X_2=1,X_1=2)}{3}$$
and 
$$\lim_{t \rightarrow \infty}\Pr(X_t=2)=\frac{2\Pr(X_2=2,X_1=1)}{3}+\frac{\Pr(X_2=1,X_1=2)}{3}+\Pr(X_2=2,X_1=2).$$}
\end{6}

Our final result connects the limiting probability distribution with the dominant eigenvectors of $Q$. A similar result can be found in \cite{WuChu17}.

\begin{2}
\label{conn}
Suppose that $\pi$ is the limiting probability distribution of an $(m-1)$th order regular Markov chain $X$ on state space $S$. Let $Q$ be the transition matrix of its reduced first order chain. Suppose that $y \in \mathbb R^N$ is a normalized nonnegative right eigenvector of $Q$ associated with the dominant eigenvalue $\lambda=1$, i.e., $Qy=y$, $y \ge 0$, and $e^Ty=1$, then, 
\be
\label{pz}
\pi=P^{(0)}y,
\ee
where $P^{(0)}$ is given in (\ref{p0}).
\end{2}
\bp
According to (\ref{komo_for0}) in Lemma \ref{komo}, we have 
$$P^{(k)}y=P^{(0)}Q^k y=P^{(0)}y.$$
From (\ref{lim_for}), we clearly have $\ds \lim_{k \rightarrow \infty}P^{(k)}=\pi e^T$. Hence 
$$\pi=P^{(0)}y,$$
i.e. (\ref{pz}) holds.
\ep

We comment that, unlike in \cite{WuChu17}, the above result holds even when $Q$ is non-regular and has multiple dominant eigenvectors associated with eigenvalues of modulus $1$. In fact, this may happen even when $\cal P$ is regular, as illustrated by the example below from \cite{Gei17}. 
\begin{6}
{\rm 
Consider a second order chain on state space $\{1,2,3,4\}$, whose transition tensor $\cal P$ has the form  
$${\cal P}(:,:,1)=\left[\ba{cccc}
1/2 & 0 & 0 & 0\\
1/2 & 0 & 1 & 0\\
0 & 1 & 0 & 1\\
0 & 0 & 0 & 0
\ea\right], ~{\cal P}(:,:,2)=\left[\ba{cccc}
0 & 0 & 1/2 & 1\\
0 & 1/2 & 0 & 0\\
1/2 & 1/2 & 0 & 0\\
1/2 & 0 & 1/2 & 0
\ea\right],$$
$${\cal P}(:,:,3)=\left[\ba{cccc}
0 & 1 & 0 & 1\\
1 & 0 & 1/2 & 0\\
0 & 0 & 1/2 & 0\\
0 & 0 & 0 & 0
\ea\right], ~{\cal P}(:,:,4)=\left[\ba{cccc}
0 & 0 & 0 & 0\\
1 & 1 & 1 & 0\\
0 & 0 & 0 & 1/2\\
0 & 0 & 0 & 1/2
\ea\right].$$
Using (\ref{pk}), we obtain ${\cal P}^{(10)} > 0$ and thus $\cal P$ is regular. Theorem \ref{conv} guarantees the existence of the limiting probability distribution $\pi$ of this chain. On the other hand, with (\ref{q}), it is easy to verify that the transition matrix $Q$ of the reduced first order chain has a double dominant eigenvalue $\lambda=1$ with an eigenspace spanned by 
$$\displaystyle y^{(1)}=\frac{1}{7}[0 ~0 ~1 ~1 ~2 ~0 ~0 ~0 ~0 ~2 ~0 ~0 ~0 ~0 ~1 ~0]^T$$ 
and 
$$\displaystyle y^{(2)}=\frac{1}{7}[0 ~2 ~0 ~0 ~0 ~0 ~2 ~0 ~1 ~0 ~0 ~1 ~1 ~0 ~0 ~0]^T.$$
Because of this, the two-phase power iteration method in \cite{WuChu17} is no longer appropriate since it does not converge to a unique answer. Nevertheless, it follows from (\ref{pz}) that 
$$\pi=P^{(0)}y^{(1)}=P^{(0)}y^{(2)}=\frac{1}{7}[2 ~2 ~2 ~1]^T.$$}
\end{6}

\section{Conclusion}
\label{concl}
\setcounter{equation}{0}

In the paper, motivated by the relevant existing works, we have developed several useful properties regarding the limiting probability distribution of a general higher order Markov chain. In particular, we have shed light on the key role the regularity of the transition tensor plays. We have also provided a way of determining the regularity via the transition matrix of the reduced first order chain. Moreover, we have connected the limiting probability distribution to some eigenvector arising on the reduced first order chain, even when this chain is non-regular.


\begin{thebibliography}{}

\bibitem{BozHaj16} H. Bozorgmanesh and M. Hajarian, Convergence of a transition probability tensor of a higher-order Markov chain to the stationary probability vector, {\it Numerical Linear Algebra with Applications} 23 (2016): 972--988.

\bibitem{ChaZha13} K.C. Chang, T. Zhang, On the uniqueness and non-uniqueness of the positive $Z$-eigenvector for transition probability tensors, {\it Journal of Mathematical Analysis and Applications} 408 (2013): 525--540.

\bibitem{CulPeaZha17} J. Culp, K. Pearson, T. Zhang, On the uniqueness of the $Z_1$-eigenvector of transition probability tensors, {\it Linear and Multilinear Algebra} 65 (2017): 891--896.

\bibitem{DinNgWei18} W. Ding, M. Ng, and Y. Wei, Fast computation of stationary joint probability distribution of sparse Markov chains, {\it Applied Numerical Mathematics} 125 (2018): 68--85.

\bibitem{FasTud20}
D. Fasino, F. Tudisco, Ergodicity coefficients for higher-order stochastic processes, {\it SIAM Journal on Mathematics of Data Science} 2 (2020): 740--769.

\bibitem{GauTudHei19}
A. Gautier, F. Tudisco, M. Hein, A unifying Perron-Frobenius theorem for nonnegative tensors via multihomogeneous maps, {\it SIAM Journal on Matrix Analysis and Applications} 40 (2019): 1206--1231.

\bibitem{Gei17} B.C. Geiger, A sufficient condition for a unique invariant distribution of a higher-order Markov chain, {\it Statistics and Probability Letters} 130 (2017): 49--56.

\bibitem{GleLimYu15} D.F. Gleich, L.H. Lim, and Y. Yu, Multilinear pagerank, {\it SIAM Journal on Matrix Analysis and Applications} 36 (2015): 1507--1541.

\bibitem{HanWanXu22} L. Han, K. Wang, and J. Xu, Higher order ergodic Markov chains and first passage times,  {\it Linear and Multilinear Algebra} 70 (2022): 6772--6779.

\bibitem{HanXu24a} L. Han and J. Xu, Ever-reaching probabilities and mean first passage times of higher order ergodic Markov chains, {\it Linear and Multilinear Algebra} 72 (2024): 59--75.

\bibitem{HanXu24b} L. Han and J. Xu, On classification of states in higher order Markov chains, {\it Linear Algebra and Its Applications} 685 (2024): 24--45.

\bibitem{HorJoh85}
R. Horn, C. Johnson, {\it Matrix Analysis}, 1st ed., Cambridge University Press, Cambridge, 1985.

\bibitem{HuQi14} S. Hu and L. Qi, Convergence of a second order Markov chain, {\it Applied Mathematics and Computation}, 241 (2014), 183--192.

\bibitem{HuaQi20} Z.H. Huang and L. Qi, Stationary probability vectors of higher-order two-dimensional symmetric transition probability tensors, {\it Asia-Pacific Journal of Operational Research} 37 (2020): 2040019, 14 pages.

\bibitem{Hun83} J.J. Hunter, {\it Mathematical Techniques of Applied Probability, Vol. 1, Discrete Time Models: Basic Theory}, Academic Press, New York, 1983.

\bibitem{Ios07} M. Iosifescu, {\it Finite Markov Processes and Their Applications}, Dover Publications, New York, 2007.

\bibitem{KemSne60} J.G. Kemeny, J.L. Snell, {\it Finite Markov Chains}, Springer-Verlag, New York, 1960.

\bibitem{KilMar11} M.E. Kilmer, C.D. Martin, Factorization strategies for third-order tensors, {\it Linear Algebra and Its Applications} 435 (2011): 641--658.

\bibitem{KolBad09} T.G. Kolda, B.W. Bader, Tensor decompositions and applications, {\it SIAM Review} 51 (2009): 455--500.

\bibitem{LiZha16} C.K. Li and S. Zhang, Stationary probability vectors of higher-order Markov chains, {\it Linear Algebra and Its Applications} 473 (2016): 114--125.

\bibitem{LiNg14} W. Li and M. Ng, On the limiting probability distribution of a transition probability tensor, {\it Linear and Multilinear Algebra} 62 (2014): 362--385.

\bibitem{MarShaLar13} C.D. Martin, R. Shafer, B. Larue, An order-$p$ tensor factorization with applications in imaging, {\it SIAM Journal on Scientific Computing} 35 (2013): 474--490.

\bibitem{QiLuo17} L. Qi and Z. Luo, {\it Tensor Analysis: Spectral Theory and Special Tensors}, SIAM, Philadelphia, 2017.

\bibitem{SmiBroGel04} A. Smilde, R. Bro, P. Geladi, {\it Multi-Way Analysis: Applications in the Chemical Sciences}, Wiley, West Sussex, UK, 2004.

\bibitem{Vla85} I. Vladimirescu, Regular homogeneous Markov chains of order two, {\it Analele Universitatii din Craiova, Seria Matematica, Fizica-Chimie} 13 (1985): 59--63 (Romanian).

\bibitem{WuChu17} S.J. Wu and M.T. Chu, Markov chains with memory, tensor formulation, and the dynamics of power iteration, {\it Applied Mathematics and Computation} 303 (2017): 226--239.

\end{thebibliography}
\end{document}